\documentclass[a4paper,11pt]{amsart}
\usepackage{amssymb}
\usepackage[all]{xy}

\theoremstyle{plain}
\newtheorem{theorem}{Theorem}[section]
\newtheorem{lemma}[theorem]{Lemma}
\newtheorem{proposition}[theorem]{Proposition}
\newtheorem{corollary}[theorem]{Corollary}
\theoremstyle{remark}
\newtheorem{remark}[theorem]{Remark}
\newtheorem{example}[theorem]{Example}

\newcommand{\dc}{\Delta}
\newcommand{\C}{\mathbb{C}}
\newcommand{\R}{\mathbb{R}}
\newcommand{\Z}{\mathbb{Z}}

\DeclareMathOperator{\sing}{Sing}

\begin{document}
\title[Structure of exceptional sets]
      {Combinatorial structure of exceptional sets
       in resolutions of singularities}
\author{D.\,A. Stepanov}
\address{The Department of Mathematical Modeling \\
         Bauman Moscow State Technical University \\
         Moscow 105005, Russia \\
         Currently: The Max-Planck-Institut f{\"u}r Mathematik \\
         Vivatsgasse 7, 53111 Bonn, Germany}
\email{dstepanov@bmstu.ru, dstepano@mpim-bonn.mpg.de}
\thanks{The research was supported by RFBR, grant no. 05-01-00353, 
        CRDF, grant no. RUM1-2692-MO-05, and the Program for the 
        Development of Scientific Potential of the Highest School, 
        no. 2.1.1.2381.}
\date{}
\begin{abstract}
The dual complex can be associated to any resolution of singularities
whose exceptional set is a divisor with simple normal crossings. It
generalizes to higher dimensions the notion of the dual graph of a
resolution of surface singularity. In this preprint we show that the
dual complex is homotopy trivial for resolutions of 3-dimensional
terminal singularities and for resolutions of Brieskorn
singularities. We also review our earlier results on resolutions of
rational and hypersurface singularities. 
\end{abstract}

\maketitle

\section{Introduction}
In this preprint we continue our study of the dual complex associated
to a resolution of singularities started in \cite{Stepanov1} and
\cite{Stepanov2}. The dual complex is defined in the
following way. Let $(X,S)$ be a germ of an algebraic variety or an
analytic space $X$, $S=\sing(X)$, and let $f\colon Y\to X$ be a good
resolution of singularities. By this we mean that the exceptional set
$Z=f^{-1}(S)$ is a divisor with simple normal
crossings on $Y$. \emph{The dual complex associated to the resolution} $f$ is
just the incidence complex $\dc(Z)$ of the divisor $Z$, i.~e., if
$Z=\sum Z_i$ is the decomposition of $Z$ into its prime components,
then $0$-simplexes (vertices) $\Delta_i$ of $\dc(Z)$ correspond to the
divisors $Z_i$, $1$-simplexes (edges) $\Delta_{ij}^{(k)}$ correspond
to the irreducible components $Z_{ij}^{(k)}$ of all intersections 
$Z_i\cap Z_j=\cup_k Z_{ij}^{(k)}$ and the edges $\Delta_{ij}^{(k)}$
join the vertices $\Delta_i$ and $\Delta_j$, 2-simplexes (triangles) 
correspond to the irreducible components of triple intersections 
$Z_i\cap Z_j\cap Z_k$ and are glued to the 1-skeleton of $\dc(Z)$ in 
a natural way, and so on. In the case $X$ to be 2-dimensional we have the usual
definition for the dual graph of resolution $f$.
\begin{example}\label{Ex:sphere}
Consider a 3-dimensional singularity
$$(\{g(x,y,z,t)=x^5+y^5+z^5+t^5+xyzt=0\},0)\subset(\C^4,0)\,.$$
Blowing up the origin produces a good resolution. The
exceptional divisor $Z$ of this blow-up is defined in the projective
space $\mathbb{P}^3$ by the homogeneous part $g_4=xyzt$ of $g$. Thus
$Z$ consists of 4 planes in general position. It follows that
$\dc(Z)$ is the surface of tetrahedron. This example can easily be
generalised to arbitrary dimension $n\geq 2$. This gives complexes
$\dc(Z)$ which are borders of standard simplexes $\Delta^{n-1}$. We see that
the dual complex associated to a resolution can be homeomorphic to the
sphere $S^{n-1}$. 
\end{example}

Note that in general $\dc(Z)$ is a triangulated topological space but
not a simplicial complex. For example, let $\dim X=2$ and $Z$ contain
2 curves meeting transversally at 2 points (see
Fig.~\ref{F:mult_edges} {\it a}). The corresponding fragment
of the dual complex is shown in Fig.~\ref{F:mult_edges} {\it b}.
\begin{center}
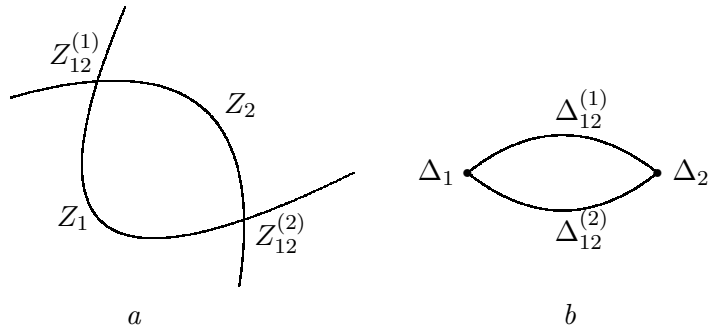
\begin{figure}[h]
\begin{picture}(10,5)
\qbezier(0.5,3.5)(4,4.5)(3.5,1)
\qbezier(2,4.7)(0,0)(5,2.5)
\put(1.1,1.8){$Z_1$} \put(3.3,3.3){$Z_2$}
\put(1,3.95){$Z_{12}^{(1)}$} \put(3.7,1.55){$Z_{12}^{(2)}$}
\put(2,0.5){\it a}
\put(6.5,2.5){\circle*{0.1}} \put(9,2.5){\circle*{0.1}}
\qbezier(6.5,2.5)(7.75,3.5)(9,2.5)
\qbezier(6.5,2.5)(7.75,1.5)(9,2.5)
\put(5.85,2.4){$\Delta_1$} \put(9.2,2.4){$\Delta_2$}
\put(7.65,3.2){$\Delta_{12}^{(1)}$} \put(7.65,1.6){$\Delta_{12}^{(2)}$}
\put(7.75,0.5){\it b}
\end{picture}
\caption{$\dc(Z)$ is not a simplicial complex.}\label{F:mult_edges}
\end{figure}
\end{center}
The complex $\dc(Z)$ is simplicial if and only if all the
intersections $Z_{i_1}\cap Z_{i_2}\cap\dots\cap Z_{i_k}$ are
irreducible. This can be achieved on some resolution of $X$. In this
case $\dc(Z)$ coincides with the topological
nerve of the covering of $Z$ by subsets $Z_i$. Also note that if
$\dim X=n$, then the condition $Z$ to have normal crossings implies
$\dim\dc(Z)\leq n-1$.

For surfaces, it is usual to consider weighted resolution graphs. The
weights assigned to vertices are the intersection numbers
$Z_i\cdot Z_i$. Also the genuses of curves $Z_i$ and the intersection
matrix $(Z_i\cdot Z_j)$ are taken into
account. We do not know what could be a generalization of this
weighting to higher dimension, so we work with the purely
combinatorial complex $\dc(Z)$. However, it is clear that $\dc(Z)$
somehow reflects the complexity of the given resolution $f$.

The complex $\dc(Z)$ can be constructed for any divisor with simple
normal crossings on some variety $Y$. If $Y$ is a K{\"a}hler manifold,
cohomologies of $\dc(Z)$ with coefficients in $\mathbb{Q}$ are interpreted as
the weight $0$ components of the mixed Hodge structure on cohomologies
of $Z$ (see~\cite{Kulikov}, Chapter 4, \S 2). In~\cite{Mumford},
D. Mumford introduced 
\emph{the compact polyhedral complex} associated to a toroidal
embedding $U\subset Y$. If $Z=Y\setminus U$ is a divisor with simple
normal crossings 
this is precisely the
dual complex $\dc(Z)$, but endowed with some additional structure. In
\cite{Mumford} it
is used in the proof of Semi-stable Reduction Theorem. An important
feature of the toroidal embedding $U\subset Y$ is that toroidal
birational morphisms
$$
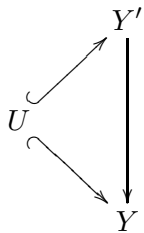
\xymatrix{
&  &Y'\ar[dd] \\
&U\ar@{^{(}->}[ru]\ar@{_{(}->}[rd] & \\
&  &Y
}
$$
are in 1-to-1 correspondence with subdivisions of the conical
polyhedral complex associated to $\dc(Z)$. We make use of this fact in
Sections~\ref{S:terminal} and \ref{S:nondeg}.

As far as we know, the first work in which the dual complex is studied
in connection to resolution of singularities (in arbitrary dimension)
is \cite{Gordon1}. There
G.\,L. Gordon considers the incidence complex of a hypersurface
$V_0\subset W$ which is the singular fiber of a map $\pi\colon W\to D$
onto the unit disc $D\subset\C$. The hypersurface $V_0$ is supposed to
have simple normal crossings. In particular, $V_0\subset W$ can be
obtained as an embedded resolution of some singular hypersurface
$H_0\subset V$, $V_0$ being the the strict transform $H'_0$ of
$H_0$ plus all exceptional divisors $Z_i$ of the embedded
resolution $W\to V$. G.\,L. Gordon shows that homologies of $\dc(V_0)$ give
some information on the monodromy around $V_0$. His article contains
also several interesting examples of the dual complex. But he
considers $\dc(V_0)$ for $V_0=H'_0+\sum Z_i$ only. According to our
definition, the dual complex associated to the resolution $H'_0$ of
$H_0$ is $\dc(\sum Z_i|_{H'_0})$.

In \cite{Shokurov} V.\,V. Shokurov studies (among other things) some
complex associated to a resolution of a log-canonical singularity
$(X,o)$. It is constructed in exactly the same way as described above
but only those prime exceptional divisors are taken into account which
have discrepancy $-1$ over $(X,o)$. This complex has a significant
property that it is uniquely determined as a topological space, i.~e., it
depends only on the singularity $(X,o)$ but not on the resolution.

The starting point of our work is the fact that homotopy type of the
dual complex $\dc(Z)$ does not depend on the choice of a resolution
$f$. Thus it is an invariant of a singularity. This was observed by
the author in \cite{Stepanov1} for isolated singularities defined over
a field of characteristic 0. The proof is an easy consequence of
Abramovich-Karu-Matsuki-W{\l}odarczyk Weak Factorization Theorem in
the Logarithmic Category (\cite{AKMW}). A. Thuillier in
\cite{Thuillier} establishes a much more general result:
\begin{theorem}\label{T:stronginv}
Let $X$ be an algebraic scheme over a perfect field $k$ and $Y$ be a
subscheme of $X$. If $f_1\colon X_1\to X$ and $f_2\colon X_2\to X$ are
two proper morphisms such that $f_{i}^{-1}(Y)$ are divisors with
simple normal crossings and $f_i$ induce isomorphisms between
$X_i\setminus f_{i}^{-1}(Y)$ and $X\setminus Y$, $i=1$, $2$, then the
topological spaces $\dc(f_{1}^{-1})$ and $\dc(f_{2}^{-1})$ are
canonically homotopy equivalent.
\end{theorem}
A. Thuillier's proof does not use the Weak Factorization
Theorem. Having in mind these results, in the sequel we sometimes simply
say that $\dc(Z)$ is the dual complex of a singularity $(X,S)$ not
indicating explicitly which good resolution is considered.

These observations motivate the following task: determine the homotopy
type of the dual complex for different classes of singularities. This
homotopy type is not always trivial as is shown by
Example~\ref{Ex:sphere}. However, it is trivial (i.~e., $\dc(Z)$ is
homotopy equivalent to a point) for many important
types of singularities. This holds, for example, if $(X,o)$ is an
isolated toric singularity (see \cite{Stepanov1}), if $(X,o)$ is a
rational surface singularity (see \cite{Artin}), or if $(X,o)$ is a
3-dimensional Gorenstein terminal singularity, defined over $\C$
(\cite{Stepanov2}). As to a general rational singularity, we can only
prove that the highest homologies of $\dc(Z)$ vanish:
$H_{n-1}(\dc(Z))=0$ (\cite{Stepanov2}). Another partial result from
\cite{Stepanov2} is that if $(X,o)$ is an isolated hypersurface
singularity over $\C$, $(Y\supset Z)\to(X\ni o)$ its good resolution,
then the fundamental group $\pi_1(\dc(Z))$ is trivial.

In this preprint we show that $\dc(Z)$ is homotopy trivial for
resolutions of all 3-dimensional terminal singularities over $\C$ (see
Theorem~\ref{T:terminal}). Applying the Varchenko-Hovanski{\u\i} 
embedded toric resolution we prove that if 
$$(Y\supset Z)\to (X\ni o)$$
is a good resolution of a non-degenerate isolated hypersurface
singularity $(X,o)$ satisfying some additional technical property, then the
complex $\dc(Z)$ does not have intermediate homologies, i.~e.,
$$H_k(\dc(Z),\Z)=0 \text{ for } 0<k<n-1,\; n=\dim X$$
(see Corollary~\ref{C:nondeg}). We describe how $\dc(Z)$ can be found for
such singularities. In particular, it follows that if $X$ is a
Brieskorn singularity, then $\dc(Z)$ is homotopy trivial
(Corollary~\ref{C:Brieskorn}).

It turns out that condition $Z$ to be the exceptional divisor of some
resolution of singularities poses strong restrictions on the homotopy
type of $\dc(Z)$. At the same time if we ask which complex $K$ can
be realized as the dual complex of some divisor $Z\subset X$ (not
necessarily contractible), then the answer is: any finite simplicial
complex. Indeed, suppose we are given a finite simplicial complex
$K$. Let $N$ be the number of its vertices and $d$ be its
dimension. First consider the divisor $H=\sum\limits_{i=1}^{N}
H_i\subset\mathbb{P}^{d+1}$ consisting of $N$ hyperplanes $H_i$ in
general position in $\mathbb{P}^{d+1}$. Its dual complex is maximal in
the sense that any $k$ vertices of $\dc(H)$, $1\leq k\leq d+1$, form a
$(k-1)$-simplex of $\dc(H)$. Thus $K$ can be identified with some
subcomplex of $\dc(H)$. Now suppose that we blow up the space
$\mathbb{P}^{d+1}$ with the center $H_{i_0}\cap H_{i_1}\cap\dots\cap
H_{i_k}$. Let $H'$ be the strict transform of $H$ under this
blow-up. Then $\dc(H')$ is obtained from $\dc(H)$ by deleting the
$k$-simplex $\Delta_{i_0 i_1\dots i_k}$ and all the simplexes to which
$\Delta_{i_0 i_1\dots i_k}$ belongs as a face. Now it is clear that
after a finite sequence of appropriate blow-ups we get a strict
transform $Z$ of $H$ such that $\dc(Z)$ is homeomorphic to $K$. We
learned this construction from \cite{Shokurov}.

The preprint is organized as follows. In Section~\ref{S:invariant} we
show that the homotopy type of $\dc(Z)$ is independent of a
resolution. The proof is reproduced from \cite{Stepanov1}; this is
done for reader's convenience. We deal mostly with isolated
singularities over $\C$, so the proof from \cite{Stepanov1} based on
the Weak Factorization Theorem is sufficient for us. In
Section~\ref{S:rh} we recall some results and sketch proofs from
\cite{Stepanov2} on the dual complex for rational and hypersurface
singularities. In Section~\ref{S:terminal} we consider resolutions of
3-dimensional terminal singularities. Section~\ref{S:nondeg} is devoted to 
non-degenerate hypersurface singularities. At the beginning of it we
recall the construction of the Varchenko-Hovanski{\u\i} embedded toric
resolution. 

The author is grateful for hospitality to the Max-Planck-Institut
f{\"u}r Mathematik in Bonn, where this preprint was written. We
especially thank to Yu.\,G. Prokhorov who pointed out a mistake in the
first version of our paper, and to K.\,A. Shramov who suggested several
improvements to the text.

\section{Invariance of the homotopy type of
  $\dc(Z)$}\label{S:invariant}
Let $(X,o)$ be a germ of an isolated singularity $o$ of an algebraic
variety or an analytic space $X$. In the algebraic case the ground
field is supposed to be of characteristic 0. In this section we prove
the following
\begin{theorem}\label{T:invariant}
(\cite{Stepanov1}) Let $f\colon Y\to X$ and $f'\colon Y'\to X$ be two
  good resolutions of $X$ with exceptional divisors $Z$ and $Z'$
  respectively. Then the dual complexes $\dc(Z)$ and $\dc(Z')$ are
  homotopy equivalent.
\end{theorem}
The proof is based on the following theorem due to
Abramovich-Karu-Matsuki-W{\l}odarczyk (see \cite{AKMW} and
\cite{Matsuki}, Theorem~5-4-1).
\begin{theorem}[Weak Factorization Theorem in the Logarithmic
    Category]\label{T:factorization}
Let $(U_{X_1},X_1)$ and $(U_{X_2},X_2)$ be complete nonsingular
toroidal embeddings over an algebraically closed field of 
characteristic zero. Let 
$$\varphi\colon(U_{X_1},X_1)--\to (U_{X_2},X_2)$$
be a birational
map which is an isomorphism over $U_{X_1}=U_{X_2}$. Then the map
$\varphi$ can be factored into a sequence of blow-ups and blow-downs
with smooth admissible and irreducible centers disjoint from
$U_{X_1}=U_{X_2}$. That is to say, there exists a sequence of
birational maps between complete nonsingular toroidal embeddings
$$(U_{X_1},X_1)=(U_{V_1},V_1)\overset{\psi_1}{-\,-\to}(U_{V_2},V_2)
\overset{\psi_2}{-\,-\to}$$
$$\dots\overset{\psi_{i-1}}{-\,-\to}(U_{V_i},V_i)
\overset{\psi_i}{-\,-\to}
(U_{V_{i+1}},V_{i+1})\overset{\psi_{i+1}}{-\,-\to}\dots$$
$$\overset{\psi_{l-2}}{-\,-\to}(U_{V_{l-1}},V_{l-1})
\overset{\psi_{l-1}}{-\,-\to}(U_{V_l},V_l)=(U_{X_2},X_2)\,,$$ 
where \par
(i) $\varphi=\psi_{l-1}\circ\psi_{l-2}\circ\dots\circ\psi_1$, \par
(ii) $\psi_i$ are isomorphisms over $U_{V_i}$, and \par
(iii) either $\psi_i$ or $\psi_{i}^{-1}$ is a morphism obtained
by blowing up a smooth irreducible center $C_i$ (or $C_{i+1}$) 
disjoint from $U_{V_i}=U_{V_{i+1}}$ and transversal to the boundary
$D_{V_i}=V_i\setminus U_{V_i}$ (or $D_{V_{i+1}}=
V_{i+1}\setminus U_{V_{i+1}}$), i. e., at each point $p\in V_i$
 (or $p\in V_{i+1}$) there exists a regular coordinate system
$\{x_1,\dots,x_n\}$ in a neighborhood $p\in U_p$ such that 
$$D_{V_i}\cap U_p\; (\text{or } D_{V_{i+1}}\cap U_p)=
\{\prod_{j\in J}x_j=0\}$$
and 
$$C_i\cap U_p\; (\text{or }C_{i+1}\cap U_p)=\{\prod_{j\in J}x_j=0\,,
\; x_{j'}=0\; \forall j'\in J'\}\,,$$
where  $J,J'\subseteq\{1,\dots,n\}$.
\end{theorem}
Here \emph{toroidal embedding} $U\subset X$ means that $U$ is an open
dense set in $X$ and $X\setminus U$ is a divisor with simple normal
crossings. Theorem~\ref{T:factorization} has also an analytic
version. To get it one needs only to replace ``birational'' with
``bimeromorphic'' etc.

In order to prove Theorem~\ref{T:invariant}, let us take the 
resolutions $(Y\setminus Z,Y)$ and $(Y'\setminus Z',Y')$ as 
toroidal embeddings and compactify $Y$ and $Y'$ to smooth
varieties (here we use the fact that the given singularity $(X,o)$
is isolated). Now Theorem~\ref{T:invariant} follows from
Theorem~\ref{T:factorization} and
\begin{lemma}
Let $\sigma\colon(X'\setminus Z',X')\to(X\setminus Z,X)$ be a
blow-up of an admissible center $C\subset Z$ in a nonsingular
toroidal embedding $(X\setminus Z,X)$, 
$X'\setminus Z'\simeq X\setminus Z$. Then the topological spaces
$\dc(Z')$ and $\dc(Z)$ have the same homotopy type.
\end{lemma}
\begin{proof}
Let $Z=\sum_{i=1}^{N}Z_i$ be the decomposition of $Z$ into its
prime components, and let $C\subset Z_i$ for $1\le i\le l$ and
$C\nsubseteq Z_i$ for $l<i\le N$. Assume that $C$ has nonempty
intersections also with $Z_{l+1},\dots,Z_r$, $l<r\le N$. There are
two possibilities. \newline
1) $\dim C=n-l$ ($n=\dim X$), i. e., $C$ coincides with one of the
irreducible components of the intersection $Z_1\cap\dots\cap Z_l$:
$C=Z_{1\dots l}^{(1)}$. Then after the blow-up the intersection of the 
proper transforms $Z'_1,\dots,Z'_l$ of the divisors $Z_1,\dots,Z_l$ 
has $J-1$ irreducible components (if $J$ is the number of components
of $Z_1\cap\dots\cap Z_l$), but all these proper transforms
intersect the exceptional divisor $F$ of the blow-up $\sigma$.
Furthermore, $F$ intersects proper transforms $Z'_{l+1},\dots,
Z'_r$ of the divisors $Z_{l+1},\dots,Z_r$. Now it is clear that
the complex $\dc(Z')$ is obtained from $\dc(Z)$ by the barycentric
subdivision of the simplex $\Delta_{1\dots l}^{(1)}$ with the center
at the point corresponding to the divisor $F$. Thus the complexes
$\dc(Z')$ and $\dc(Z)$ are even homeomorphic. \newline
2) $\dim C<n-l$, let $C\subset Z_{1\dots l}^{1}$. In this case
divisors $Z_{i_1},\dots,Z_{i_s}$ have nonempty intersection if
and only if their proper transforms $Z'_{i_1},\dots,Z'_{i_s}$ have
nonempty intersection. Therefore the complex $\dc(Z')$ is
obtained from the complex $\dc(Z)$ in the following way. Add to $\dc(Z)$ 
a new vertex corresponding to the exceptional divisor $F$ of the
blow-up $\sigma$ and construct cones with vertex at $F$ over all
the maximal cells $\Delta_{i_1\dots i_s}^{(j)}$ of the complex $\dc(Z)$
possessing the property
$$Z_{i_1\dots i_s}^{(j)}\cap C\ne \varnothing\,.$$
Note that the simplex $\Delta_{F,1\dots l}$ corresponding to
the intersection $F\cap Z_{1\dots l}^{(1)}$ is regarded as
a common simplex for all constructed cones. Now we can define the
homotopy equivalence between $\dc(Z')$ and $\dc(Z)$ as a contraction
of the constructed cones: it sends the vertex $F$ of the complex 
$\dc(Z')$ to any of the vertices $Z_1,\dots,Z_l$ of the cell
$\Delta_{1\dots l}^{(1)}$ of the complex $\dc(Z)$ and it is identity
on other vertices of $\dc(Z')$ ($\dc(Z)$). Then the induced
simplicial map is our homotopy equivalence. 
\end{proof} 

Figure~\ref{F:invariant} illustrates part 2) of the proof. Here we
suppose that $\dim X=3$ and $Z$ consists of $4$ prime components. We
denote the corresponding vertices of $\dc(Z)$ by the same letters
$Z_i$. Let $\dc(Z)$ be as shown in Figure~\ref{F:invariant}, {\it
  a}. Then let us blow up a smooth irreducible curve $C\subset Z_1$ which
intersect transversally the curves $Z_1\cap Z_i$, $i=2$, $3$, $4$. For
the simplicity of drawing we assume that $C$ intersects every $Z_1\cap
Z_i$ at a single point. Then the obtained $\dc(Z')$ is shown in
Figure~\ref{F:invariant}, {\it b}. Here all the triangles belong to
$\dc(Z')$ together with their interiors.
\begin{center}
\begin{figure}[htb]
\begin{picture}(10,5)
\put(2.5,2.5){\circle*{0.1}} \put(2,1.5){\circle*{0.1}}
\put(1.5,3){\circle*{0.1}} \put(3.5,3){\circle*{0.1}}
\put(2.5,2.5){\line(-1,-2){0.5}} \put(2.5,2.5){\line(-2,1){1}}
\put(2.5,2.5){\line(2,1){1}}
\put(2.75,2.25){$Z_1$} \put(1.4,1.4){$Z_2$}
\put(1,3.2){$Z_3$} \put(3.6,3){$Z_4$}
\put(2.4,0.5){\it a}

\put(4,2){\vector(1,0){1}}

\put(6.5,2){\circle*{0.1}} \put(6,1){\circle*{0.1}}
\put(5.5,2.5){\circle*{0.1}} \put(7.5,2.5){\circle*{0.1}}
\put(6.5,4){\circle*{0.1}}
\put(6.5,2){\line(-1,-2){0.5}} \put(6.5,2){\line(-2,1){1}}
\put(6.5,2){\line(2,1){1}}
\put(6.5,4){\line(0,-1){2}} \put(6.5,4){\line(-1,-6){0.5}}
\put(6.5,4){\line(-2,-3){1}} \put(6.5,4){\line(2,-3){1}}
\put(6.75,1.75){$Z'_1$} \put(5.4,0.9){$Z'_2$}
\put(5,2.8){$Z'_3$} \put(7.6,2.5){$Z'_4$}
\put(6.5,4.1){$F$}
\put(6.4,0.5){\it b}
\end{picture}
\caption{Transformation of $\dc(Z)$ into $\dc(Z')$}\label{F:invariant}
\end{figure}
\end{center}

\section{Dual complex for rational and hypersurface
  singularities}\label{S:rh}
In this section and in the rest of the paper we consider only
varieties (and analytic spaces) over $\C$. Also when we speak about
the dual complex $\dc(Z)$ associated to a resolution of a given
singularity $(X,o)$, we shall always assume that $\dc(Z)$ is a
simplicial complex. This can be achieved on some resolution of $X$ and
we already established the invariance of the homotopy type of
$\dc(Z)$ in section~\ref{S:invariant}.

\subsection{Rational singularities}
First recall that an algebraic variety (or an analytic space) $X$ has 
\emph{rational singularities} if $X$ is normal and for any resolution 
$f\colon Y\to X$ all the sheaves $R^i f_*\mathcal{O}_{Y}$ vanish, 
$i>0$.

It is well known that the exceptional divisor in a resolution of a
rational surface singularity is a tree of rational curves. This
follows, e.~g., from M. Artin's considerations in \cite{Artin}. Thus
in the surface case the dual graph for a rational singularity is
homotopy trivial. In \cite{Stepanov2} we partially generalized this
statement to higher dimensions and proved the following result.
\begin{theorem}\label{T:rational}
Let $o\in X$ be an isolated rational singularity of a variety 
(or an analytic space) $X$ of dimension $n\geq 2$, and let 
$f\colon Y\to X$ be a good resolution with the exceptional divisor 
$Z$. Then the highest homologies of the complex $\dc(Z)$ 
vanish:
$$H_{n-1}(\dc(Z),\mathbb{Z})=0\,.$$
\end{theorem}
In the first step of the \emph{proof} we follow M. Artin's argument
from \cite{Artin}. Let $Z=\sum\limits_{i=1}^{N} Z_i$ be the decomposition of the 
divisor $Z$ into its prime components $Z_i$. We can assume that $X$ 
is projective (since the given singularity is isolated) and $f$
is obtained by a sequence of smooth blow-ups (Hironaka's resolution 
\cite{Hironaka}). Thus all $Z_i$ and $Y$ are K{\"a}hler manifolds.

The sheaves $R^i f_*\mathcal{O}_{Y}$  are concentrated at the 
point $o$. Via Grothendieck's theorem on formal functions 
(see \cite{Grothendieck}, (4.2.1), and \cite{Grauert}, Ch. 4,
Theorem~4.5 for the analytic case) the completion of the stalk of 
the sheaf $R^i f_*\mathcal{O}_{Y}$ at the point $o$ is
\begin{equation}\label{E:projlim}
\varprojlim_{(r)\to(\infty)} H^i(Z,\mathcal{O}_{Z_{(r)}})\,,
\end{equation}
where $(r)=(r_1,\dots,r_N)$ and $Z_{(r)}=\sum\limits_{i=1}^{N}
r_iZ_i$. If $(r)\geq (s)$, i.~e., $r_i\geq s_i$ $\forall i$, there 
is a natural surjective map $g$ of sheaves on $Z$: 
$$g\colon\mathcal{O}_{Z_{(r)}}\to\mathcal{O}_{Z_{(s)}}\,.$$ 
Since dimension of $Z$ is $n-1$, the map $g$ induces a surjective 
map of cohomologies
$$H^{n-1}(Z,\mathcal{O}_{Z_{(r)}})\to
H^{n-1}(Z,\mathcal{O}_{Z_{(s)}})\,.$$
Recall that the given singularity $o\in X$ is rational, and thus the
projective limit \eqref{E:projlim} is $0$. Therefore the 
cohomology group $H^{n-1}(Z,\mathcal{O}_Z)$ vanishes too (because
the projective system in \eqref{E:projlim} is surjective). For
surfaces, $n=2$, thus $H^1(Z,\mathcal{O}_Z)=0$. It easy follows, e.~g.,
from exponential exact sequence, that $H^1(Z,\Z)=0$ and this implies
$$H^1(\dc(Z),\Z)=H_1(\dc(Z),\Z)=0\,.$$
For an arbitrary $n\geq 2$ we need a more sophisticated argument.
\begin{lemma}\label{L:cohomology}
Let $Z=\sum Z_i$ be a reduced divisor with simple normal crossings
on a compact K{\"a}hler manifold $Y$, $\dim Y=n\geq 2$, and assume that
$H^k(Z,\mathcal{O}_Z)=0$ for some $k$, $1\leq k\leq n-1$. Then the 
$k$-th cohomologies with coefficients in $\mathbb{C}$ of the complex
$\dc(Z)$ vanish too:
$$H^k(\dc(Z),\mathbb{C})=0\,.$$
\end{lemma}
For the proof see \cite{Stepanov2}. The idea is to
introduce a kind of Mayer-Vietoris spectral sequence for $Z=\cup Z_i$
and to show that it degenerates in $E_2$. Compare also \cite{Gordon2}
and \cite{Kulikov} 
where very closed results are stated.

Lemma~\ref{L:cohomology} implies $H^{n-1}(\dc(Z),\C)=0$; hence 
$H_{n-1}(\dc(Z),\Z)=0$. This completes the proof of
Theorem~\ref{T:rational}.

So, the highest homologies of $\dc(Z)$ vanish for rational
singularities. In general we do not know anything about intermediate
homologies. If we could prove that they vanish together with the
fundamental group $\pi_1(\dc(Z))$, then it would follow that $\dc(Z)$
is homotopy equivalent to a point. We can prove this only in the
partial cases of 3-dimensional hypersurface singularities (see the
next subsection) and for some hypersurfaces which are non-degenerate
in Varchenko-Hovanski{\u\i} sense (see
section~\ref{S:nondeg}). Moreover, we do not know any example of an
isolated $n$-dimensional (not necessarily rational) singularity such
that $H_k(\dc(Z))\ne 0$ for $0<k<n-1$.

\subsection{Hypersurface singularities}
If $(X,0)\subset(\C^{n+1},0)$ is a germ of an isolated hypersurface
singularity, $n\geq 3$, then its link (the intersection of $X$ with a
sufficiently small sphere around $0\in\C^{n+1}$) is simply connected
(J. Milnor \cite{Milnor}). This implies the following 
result (\cite{Stepanov2}):
\begin{theorem}\label{T:hypersurface}
Let $o\in X$ be an isolated hypersurface singularity of an
algebraic variety (or an analytic space) $X$ of dimension at least 
$3$ defined over the field $\mathbb{C}$ of complex numbers. If 
$f\colon Y\to X$ is a good resolution of $o\in X$, $Z$ its 
exceptional divisor, then the fundamental group of $\dc(Z)$ is trivial:
$$\pi(\dc(Z))=0\,.$$
\end{theorem}
For the \emph{proof}, it suffices to notice that there is a continuous
map with connected fibers from the link $M$ of $(X,o)$ onto the
exceptional divisor $Z$ (see \cite{ACampo}). This gives
$\pi_1(Z)=0$. But we also can construct a map $\psi$ from $Z$ to $\dc(Z)$
with the same properties. It is defined as follows. 

Let us take a triangulation $\Sigma'$ of $Z$ such that all 
the intersections $Z_{i_0\dots i_p}$ are subcomplexes. Next we make 
the barycentric subdivision $\Sigma$ of $\Sigma'$ and the 
barycentric subdivision of the complex $\dc(Z)$. Now let $v$ be a 
vertex of $\Sigma$ belonging to the subcomplex $Z_{i_0\dots i_p}$ 
but not to any smaller subcomplex $Z_{i_0\dots i_p i_{p+1}}$:
$$v\in Z_{i_0\dots i_p}\,,\quad v\notin Z_{i_0\dots i_p i_{p+1}}
\:\forall\,i_{p+1}\,.$$ 
Then let 
$$\psi(v)=\text{ the center of the simplex }
\Delta_{i_0\dots i_p}\,.$$
This determines the map $\psi$ completely
as a simplicial map (depending on the triangulation $\Sigma'$).

It easily follows that $\psi$ is continuous, surjective, and has
connected fibers. Thus $\pi_1(\dc(Z))=0$.

Combining Theorems~\ref{T:rational} and \ref{T:hypersurface} we obtain
\begin{corollary}\label{C:dim3}
Let $o\in X$ be an isolated rational hypersurface singularity
of dimension $3$. If $f\colon Y\to X$ is a good resolution with the
exceptional divisor $Z$, then the dual complex $\dc(Z)$ associated
to the resolution $f$ has the homotopy type of a point.
\end{corollary}
\begin{proof}
We know from Theorems~\ref{T:rational} and \ref{T:hypersurface}
that $\dc(Z)$ is simply connected and $H_2(\dc(Z),\Z)=0$. 
Since $\dim X=3$, we have $\dim(\dc(Z))\leq 2$. Now 
Corollary~\ref{C:dim3} follows from 
the Inverse Hurevicz and Whitehead Theorems.
\end{proof}

\section{Dual complex for 3-dimensional terminal singularities}
\label{S:terminal}
Terminal singularities arose in the framework of Mori theory as
singularities which can appear on minimal models of algebraic
varieties of dimension $\geq 3$. In dimension $3$ terminal
singularities are completely classified up to an analytic equivalence
by M. Reid, D. Morrison, G. Stevens, S. Mori and
N. Shepherd-Barron. Classification looks as follows. Gorenstein (or
\emph{index} $1$) terminal singularities are exactly isolated compound
Du Val (cDV) points. A \emph{cDV-point} is a germ $(X,o)$ of
singularity analytically isomorphic to the germ
$$(\{f(x,y,z)+tg(x,y,z,t)=0\},0)\subset(\C^4,0)\,,$$
where $f$ is one of the following Klein polynomials
$$x^2+y^2+z^{n+1},\;n\geq 1\,,\: x^2+y^2z+z^{n-1},\;n\geq 4\,,$$
$$x^2+y^3+z^4,\:x^2+y^3+yz^3,\:x^2+y^3+z^5\,.$$
In other words, a cDV-point is a germ of singularity such that its
general hyperplane section is a Du Val point. Non-Gorenstein (index
$\geq 2$) terminal singularities are quotients of isolated cDV-points
by some cyclic group actions, see \cite{RY} or \cite{Mori} for a
precise statement.

\begin{theorem}\label{T:terminal}
Let $f\colon(Y,Z)\to(X,o)$ be a good resolution of a three-dimensional
terminal singularity $(X,o)$. Then the dual complex $\dc(Z)$ of $f$ is
homotopy trivial.
\end{theorem}
\begin{proof}
First suppose that the singularity $(X,o)$ is Gorenstein. Then,
according to the classification, it is an isolated hypersurface
singularity. On the other hand, all terminal (and, moreover,
canonical) singularities are rational (R. Elkik \cite{Elkik}). Thus
Corollary~\ref{C:dim3} applies and we get the homotopy triviality of
$\dc(Z)$.

Now assume that $(X,o)$ is a non-Gorenstein terminal singularity of
index $m$. Let $(V,o')\to(X,o)$ be its Gorenstein cover, $X=V/\Z_m$,
and consider the diagram
$$
\xymatrix{
&**[l]Z'\subset W\ar[r]\ar[d]_{g'} &**[r]\widetilde{X}\ar[d]^g 
\supset \widetilde{Z} \\
&**[l]o' \in V\ar[r] &**[r]X \ni o
}
$$
where $g'\colon W\to V$ is an equivariant Hironaka resolution of $V$
(see, e.~g., \cite{Hauser}), $\widetilde{X}=W/\Z_m$, horizontal arrows
stand for natural projections, and $g\colon\widetilde{X}\to X$ is the
induced birational morphism. By $Z'$ and $\widetilde{Z}$ we denote the
exceptional divisors of $g'$ and $g$ respectively. 

Since $\widetilde{X}$ has only cyclic quotient singularities,
$\widetilde{X}\setminus\widetilde{Z}\subset\widetilde{X}$ is a
toroidal embedding in the sense of \cite{Mumford}, Chapter II, \S1,
Definition~1. It is clear that the definition of the dual complex
applies also to the partial resolution $g$. It is also possible to
interpret $\dc(\widetilde{Z})$ as the underlying CW-complex of
Mumford's compact polyhedral complex of the toroidal embedding
$\widetilde{X}\setminus\widetilde{Z}\subset\widetilde{X}$ (see
\cite{Mumford}, pp. 69--71). On the other hand, the action of $\Z_m$
on $W$ naturally induces an action of $\Z_m$ on $\dc(Z')$, so that
$\dc(\widetilde{Z})=\dc(Z')/\Z_m$. 

We already know that $\dc(Z')$ is homotopy trivial. It is a
topological fact that the quotient of a finite homotopy trivial
CW-complex by a finite group is again homotopy trivial
(\cite{tomDieck}, p. 222, Theorem~6.15). Thus we get homotopy
triviality of $\dc(\widetilde{Z})$.

Let $\dc'(\widetilde{Z})$ be the conical polyhedral complex associated
to $\dc(\widetilde{Z})$. Any subdivision of $\dc'(\widetilde{Z})$
gives rise to a new toroidal embedding $Y\setminus Z\subset Y$ and a
birational toroidal morphism $(Y\setminus
Z,Y)\to(\widetilde{X}\setminus\widetilde{Z},\widetilde{X})$
(\cite{Mumford}, Chapter II, \S2, Theorem~$6^*$). The corresponding dual
complex $\dc(Z)$ is a subdivision of $\dc(\widetilde{Z})$. In
particular, we can take a subdivision of $\dc'(\widetilde{Z})$ which
gives a good resolution $(Y,Z)$ of $\widetilde{X}$. Since $\dc(Z)$ is
just a subdivision of $\dc(\widetilde{Z})$, it is also homotopy
trivial.
\end{proof}

\section{Non-degenerate hypersurface singularities}\label{S:nondeg}
\subsection{Varchenko-Hovanski{\u\i} embedded toric resolution}
We shall connect the dual complex of a non-degenerate hypersurface
singularity to its Newton diagram. For this we need to recall  the
definition of non-degeneracy and the construction of
Varchenko-Hovanski{\u\i} embedded toric resolution. We mainly follow
\cite{Varchenko}. 

Let $f\in\C\{x_0,\dots,x_n\}$ be a convergent power
series,
$$f(x)=\sum_{m\in\Z_{\geq 0}^{n+1}} a_m x^m\,,$$
where $x^m=x_{0}^{m_0}\dots x_{n}^{m_n}$. \emph{The Newton polyhedron}
$\Gamma_+(f)$ of $f$ is the convex hull of the set
$\bigcup\limits_{a_m\ne 0}(m+\R_{\geq 0}^{n+1})$ in $\R^{n+1}$. The union
of all compact faces of $\Gamma_+(f)$ is \emph{the Newton diagram}
$\Gamma(f)$ of $f$. If $\gamma$ is a face of $\Gamma(f)$, then we
denote by $f_\gamma$ the polynomial consisting of all terms $a_m x^m$
of $f$ such that $m\in\gamma$, i.~e.,
$$f_\gamma(x)=\sum_{m\in\gamma} a_m x^m\,.$$
A series $f\in\C\{x_0,\dots,x_n\}$ is called \emph{non-degenerate} if
for any face $\gamma$ of $\Gamma(f)$ the polynomial
$f_\gamma(x_0,\dots,x_n)$ defines a non-singular hypersurface in
$(\C^*)^{n+1}$.

For notation and general background in toric geometry see, e.~g.,
\cite{Danilov} or \cite{Fulton}. Suppose that 
$$(X,0)=(\{f(x_0,\dots,x_n)=0\},0)\subset(\C^{n+1},0)$$
is an isolated hypersurface singularity given by a non-degenerate
series $f$. Let us consider the ambient space $\C^{n+1}$ as a toric
variety $V(\sigma_+)$, where $\sigma_+$ is the non-negative cone
$\R_{\geq 0}^{n+1}$ in $\R^{n+1}$. A.\,N. Varchenko in
\cite{Varchenko} constructs a subdivision $\Sigma$ of the cone
$\sigma_+$ such that the corresponding toric morphism $\pi\colon
V(\Sigma)\to V(\sigma_+)$ is an embedded resolution of $(X,0)$. Here this
means that (i) $V(\Sigma)$ and the strict transform $Y$ of $X$ are
smooth, and (ii) the union of $Y$ and the exceptional set $Z$ of $\pi$
is a divisor with simple normal crossings. 

Denote by $W$ the space $\R^{n+1}$ where the cone $\sigma_+$ lies and
by $W^*$ its dual. We shall consider the Newton diagram $\Gamma(f)$ as
a subset of $W^*$. For any $w\in\sigma_+$ we can associate a number
$\mu(w)=\min\limits_{m\in\Gamma_+(f)}\langle w,m\rangle$, where
$\langle\cdot,\cdot\rangle$ stands for the pairing between $W$ and
$W^*$, and a face 
$$\gamma(w)=\{m\in\Gamma_+(f)|\langle
w,m\rangle=\mu(w)\}$$ of the Newton polyhedron $\Gamma_+(f)\,.$ 
Two vectors $w^1$ and $w^2\in \sigma_+$ are called \emph{equivalent with
respect to} $\Gamma_+(f)$ if they cut the same face on
$\Gamma_+(f)$. We shall write this equivalence relation as $w^1\sim_f
w^2$, so that
$$w^1\sim_f w^2 \Longleftrightarrow \gamma(w^1)=\gamma(w^2)\,.$$
It is not difficult to verify that closures of equivalence classes of
$\sim_f$ are rational polyhedral cones and these cones posses all the
properties necessary to form a fan. Denote this fan by $\Sigma'$. It
will be called \emph{the first Varchenko subdivision of}
$\sigma_+$. To get $\Sigma$, subdivide $\Sigma'$ so that all cones of
$\Sigma$ give non-singular affine toric varieties. This is equivalent
to saying that every cone $\sigma\in\Sigma$ is simplicial and its
skeleton (the set of primitive vectors of $\Z^{n+1}\subset W$ along
the edges of $\sigma$) is a part of a basis for $\Z^{n+1}$. The toric
variety $V(\Sigma)$ is smooth by the construction. The rest of needed
properties of the birational morphism $\pi\colon V(\Sigma)\to\C^{n+1}$ 
follow from
\begin{lemma}\label{L:Varchenko}
(i) The strict transform $Y$ of $X$ by the morphism $\pi$ is
  smooth;
(ii) all the toric strata $Z_\sigma$ of $V(\Sigma)$ corresponding to
  the cones $\sigma$ of $\Sigma$ are transversal to $Y$ in some
  neighborhood of $\pi^{-1}(0)$.
\end{lemma}
\begin{proof}
(i) Let us take some of the affine pieces of $V(\Sigma)$, say
$V(\sigma)\simeq\C^{n+1}$, corresponding to the $n+1$-dimensional cone
  $\sigma$ of $\Sigma$. Let $w^0$, $w^1,\dots,w^n$, be the skeleton of
  $\sigma$,
$$w^0=(w_{0}^{0},w_{1}^{0},\dots,w_{n}^{0}),\dots, 
w^n=(w_{0}^{n},w_{1}^{n},\dots,w_{n}^{n})\,.$$
Then the map $\pi$ restricted to $V(\sigma)$ is given by the formulae
\begin{align}
x_0 &=y_{0}^{w_{0}^{0}}y_{1}^{w_{0}^{1}}\dots y_{n}^{w_{0}^{n}}\,,
\notag \\
x_1 &=y_{0}^{w_{1}^{0}}y_{1}^{w_{1}^{1}}\dots y_{n}^{w_{1}^{n}}\,,
\label{E:tmap} \\
    &\dots \dots \notag \\
x_n &=y_{0}^{w_{n}^{0}}y_{1}^{w_{n}^{1}}\dots y_{n}^{w_{n}^{n}}\,, \notag
\end{align}
where $y_i$ are the coordinates on $V(\sigma)$. The full transform of
$X$ is given in $V(\sigma)$ by the equation
$$y_{0}^{\mu(w^0)}y_{1}^{\mu(w^1)}\dots
  y_{n}^{\mu(w^n)}f'(y_0,y_1,\dots,y_n)=0\,,\quad f'(0,0,\dots,0)\ne
  0\,,$$
and the strict transform $Y$ is $\{f'(y_0,\dots,y_n)=0\}$.

Suppose that $Y$ is singular in some point
$Q=(y_{0}^{0},\dots,y_{n}^{0})$. We assume that
$y_{0}^{0},\dots,y_{k}^{0}\ne 0$, $y_{k+1}^{0}=\dots=y_{n}^{0}=0$,
$0\leq k\leq n$. We can write
$$f'(y_0,\dots,y_k)=g(y_0,\dots,y_k)+y_{k+1}(\dots)+\dots+y_n(\dots)\,.$$
Here $g$ comes from those monomials $x^m$ of $f$ whose degrees
$\langle w^j,m\rangle$ with respect to weight-vectors
$w^{k+1},\dots,w^n$ are exactly $\mu(w^{k+1}),\dots,\mu(w^n)$. In
other words, $g=f'_{\gamma}$ (strict transform of $f_\gamma$) for the
face
$$\gamma=\gamma(w^{k+1})\cap\gamma(w^{k+2})\cap\dots\cap\gamma(w^n)\,.$$

At the same time the point $Q$ must lie on the exceptional divisor of
$\pi$. Thus some of $\mu(w^{k+1}),\dots,\mu(w^n)$ are strictly
positive. It follows that the face $\gamma$ is compact and hence
$f'_{\gamma}$ is a polynomial.

The hypersurface $\{f'_{\gamma}(y_0,\dots,y_k)=0\}$ is singular at
the point $Q$ and, moreover, in every point
$(y_{0}^{0},\dots,y_{k}^{0},y_{k+1},\dots,y_n)$ for arbitrary numbers
$y_{k+1},\dots,y_n\in\C$. Take, for instance, a point
$$Q'=(y_{0}^{0},\dots,y_{k}^{0},y'_{k+1},\dots,y'_n)$$ 
for some
$y'_{k+1},\dots,y'_n\ne 0$. But the morphism $\pi$ is a local
isomorphism at $Q'$, thus the hypersurface $\{f_\gamma=0\}$ is
singular at the point $P=\pi(Q)\in(\C^*)^{n+1}$, but this contradicts
the non-degeneracy.

(ii) In the notation of the case (i), let us show that $Y$ is
transversal to the toric stratum $L=\{y_{k+1}=\dots=y_n=0\}$. Take
some point $Q\in L\cap Y$, $\pi(Q)=0$. It is sufficient to prove that
one of the partial derivatives
$$\frac{\partial f'}{\partial y_i}(Q)\ne 0\,,\quad i=1,\dots,k\,.$$
Let $Q=(y_{0}^{0},\dots,y_{l}^{0},0\dots,0)$, $0\leq l\leq k$. Now the
same argument as in the case (i) shows that non-transversality of $L$
and $Y$ at $Q$ would contradict non-degeneracy.
\end{proof}
\begin{remark}\label{R:Varchenko}
The assertions of Lemma~\ref{L:Varchenko} essentially hold and can be
proved in the same manner also in the case when $\sigma$ is an
$(n+1)$-dimensional simplicial cone contained in one of the cones of
the fan $\Sigma'$, and we take $\pi$ to be the birational morphism
$V(\sigma)\to\C^{n+1}$. Now we do not assume the skeleton of $\sigma$
to be a part of a basis, so that $V(\sigma)$ is a quotient of $\C^{n+1}$ by
some abelian group $G$. We have to replace ``smooth'' of (i) by
``quasismooth'' indicating that the cover $\widetilde{Y}$ of $Y$ is
smooth in $\C^{n+1}$. The transversality of (ii) also must be
understood as transversality of $\widetilde{Y}$ and toric strata in
the covering $\C^{n+1}$.
\end{remark}

\subsection{Dual complex for non-degenerate singularities}
Naturally, the dual complex $\dc(Z)$ associated to a resolution of a
non-degenerate isolated hypersurface singularity
$$(X,o)=(\{f=0\},0)\subset(\C^{n+1},0)$$ 
is connected to its Newton
diagram $\Gamma(f)$. However, one should be careful applying the
embedded toric resolution to the calculation of $\dc(Z)$ because
Varchenko-Hovanski{\u\i} resolution is not a resolution in rigorous
sense. When one says that $\pi\colon (Y,Z)\to(X,o)$ is a resolution of
singularity $(X,o)$, it is usually meant that $Y\setminus Z\simeq
X\setminus\{o\}$. But Varchenko-Hovanski{\u\i} resolution can involve
blow-ups with centers different from the singular point $o$. Indeed, look at the
example  
$$(V,0)=(\{x^4+y^4+xz+yz=0\},0)\subset(\C^3,0)\,.$$
This is a non-degenerate isolated surface singularity. Its Newton
polyhedron is shown in Fig.~\ref{F:badNdiag}.
\begin{center}
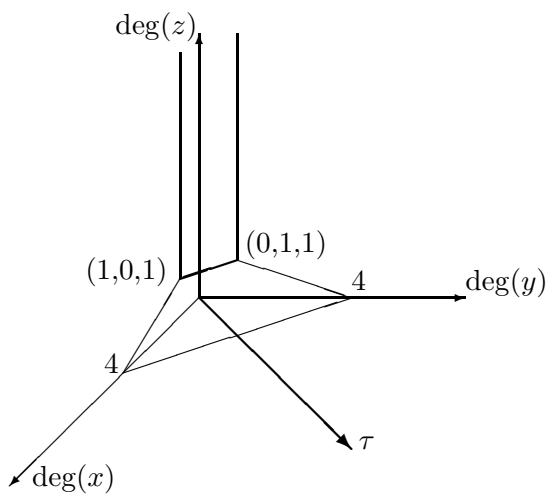
\begin{figure}[h]
\begin{picture}(7,7)
\put(3,3){\vector(-1,-1){2.5}} \put(0.8,0.5){$\deg(x)$}
\put(3,3){\vector(1,0){3.5}} \put(6.5,3.1){$\deg(y)$}
\put(3,3){\vector(0,1){3.5}} \put(1.9,6.5){$\deg(z)$}
\put(2,2){\line(3,1){3}} \put(2,2){\line(3,5){0.75}}
\put(3.5,3.5){\line(3,-1){1.5}}
\thicklines
\put(3.5,3.5){\line(-3,-1){0.75}} 
\put(2.75,3.25){\line(0,1){3}} \put(3.5,3.5){\line(0,1){3}}
\put(3,3){\vector(1,-1){2}}
\thinlines
\put(1.75,2){4} \put(5,3.1){4}
\put(1.5,3.25){(1,0,1)} \put(3.6,3.6){(0,1,1)}
\put(5.1,1){$\tau$}
\end{picture}
\caption{Newton diagram for $x^4+y^4+xz+yz$.}\label{F:badNdiag}
\end{figure}
\end{center}
The corresponding first Varchenko subdivision $\Sigma'$ contains the
ray $\tau=\langle(1,1,0)\rangle$ ``orthogonal'' to the face bordered
by bold lines in Fig.~\ref{F:badNdiag}. The cone $\tau$ gives rise to
some exceptional divisor $Z_\tau$ in Varchenko-Hovanski{\u\i}
resolution of $(V,0)$; the center of $Z_\tau$ in $\C^3$ is the
$z$-axis which is also contained in $V$.

Let us recall from \cite{Stepanov1} how $\dc(Z)$ looks when $Z$ is
the exceptional divisor of some toric birational morphism $\pi\colon
V(\Sigma)\to V(\sigma_+)=\C^{n+1}$ corresponding to a subdivision
$\Sigma$ of $\sigma_+$ such that all cones of $\Sigma$ are
simplicial. Prime exceptional divisors of
$\pi$ are in 1-to-1 correspondence with 1-dimensional cones (rays) of
$\Sigma$. Two such divisors intersect iff corresponding rays are faces
of some 2-dimensional cone of $\Sigma$. Add to $Z$ the divisor
$T=\sum\limits_{i=0}^{n} T_i$, where $T_i$ correspond to the ray
$\langle e_i=(0,\dots,\underset{i}{1},\dots,0)\rangle$. If we take
a hyperplane $H$ in $\R^{n+1}$ such that it intersects all the rays
$\langle e_i\rangle$ in a point different from $0$, then we get a
compact polyhedron $K=H\cap \sigma_+$. The fan $\Sigma$ determines
some triangulation of $K$. It is clear that $K$ with this triangulation
is exactly $\dc(Z+T)$. To obtain $\dc(Z)$, we have only to throw away
the vertices $T_i$ and all incidental to them cones.

One may object that earlier the complex $\dc(Z)$ was introduced only
for divisors with simple normal crossings, but here even the ambient
space $V(\Sigma)$ can be singular. But since all the cones of $\Sigma$
are simplicial, we can understand simple normal crossings in the
``orbifold sense'' as in Remark~\ref{R:Varchenko}, or consider
$\dc(Z)$ as the underlying complex of the polyhedral complex
associated to the toroidal (here simply toric) embedding
$V(\Sigma)\setminus Z\subset V(\Sigma)$.

Now let us come back to the singularity $(X,o)$. Assume that
\begin{description}
\item[(R)] the first Varchenko subdivision $\Sigma'$ for $(X,o)$ does
  not contain any rays on the border of cone $\sigma_+$ with the exception
  of $\langle e_0\rangle$, $\dots,\langle e_n\rangle$.
\end{description}
We shall refer to this as to the property (R). This guarantees that
the second subdivision $\Sigma$ can also be chosen with the property
(R), so that the obtained Varchenko-Hovanski{\u\i} resolution is a
resolution in rigorous sense.

Denote by $\sigma_1,\dots,\sigma_N$ all the $(n+1)$-dimensional cones
of $\Sigma'$ and let us consider one more subdivision $\Sigma''$ of
$\Sigma'$ satisfying property (R) and such that all its cones are
simplicial. We do not demand $V(\Sigma'')$ to be smooth, but
Lemma~\ref{L:Varchenko} applies to $\Sigma''$ (see
Remark~\ref{R:Varchenko}).
\begin{proposition}\label{P:intersection}
Let $\tau$ be a cone of $\Sigma''$ such that the corresponding toric
stratum $Z_\tau$ of $V(\Sigma'')$ is exceptional. Denote by $Y''$ the
strict transform of $X$ in $V(\Sigma'')$. Then $Y''\cap
Z_\tau\ne\varnothing$ iff $\tau$ does not contain any point from
interior of some $\sigma_i$, $i=1,\dots,N$. If $\dim Z_\tau\geq 2$,
then the intersection $Y''\cap Z_\tau$ is irreducible.
\end{proposition}
\begin{proof}
The fact that if $\tau$ contains a point in the interior of some
$\sigma_i$, then $Y''\cap Z_\tau=\varnothing$ follows from the
construction of $\Sigma'$. Thus suppose that $\tau$ is contained in an
$n$-dimensional cone $\sigma'$ of $\Sigma'$. Also let
$\tau\subset\sigma$, where $\sigma$ is one of the
$(n+1)$-dimensional cones of $\Sigma''$. In the affine piece
$V(\sigma)$ the morphism $\pi''\colon V(\Sigma'')\to
V(\sigma_+)\simeq\C^{n+1}$ is given by formulae~\eqref{E:tmap}, where
$y_i$ are now the coordinates on $\C^{n+1}$ which covers $V(\sigma)$
(see Remark~\ref{R:Varchenko}). Now $Z_\tau$ is defined, say, by
$y_0=\dots= y_k=0$, and the intersection $Z_\tau\cap Y''$ by the
system
$$
\begin{cases}
y_i=0,\; 0\leq i\leq k\,, \\
g(y_{k+1},\dots,y_n)=f'(0,\dots,0,y_{k+1},\dots,y_n)=0\,,
\end{cases}
$$
where $f'$ is the equation (of the cover) of the strict transform of
$X$.

Let $w^0,\dots,w^k$ be the skeleton of $\tau$. The polynomial $g$
comes from those monomials of $f$ which belong to the face
$\gamma=\gamma(w^0)\cap\dots\cap\gamma(w^k)$. But all $w_i$ belong to
an $n$-dimensional cone of $\Sigma'$, thus $\gamma$ contains at least
$2$ vertices of $\Gamma(f)$. Therefore $Z_\tau\cap Y''$ is indeed
non-empty. At the same time, $Z_\tau$ is transversal to $Y''$ in the sense of
Remark~\ref{R:Varchenko}, and, moreover, $Y''$ is transversal to all
the exceptional strata. This shows that $Z_\tau\cap Y''$ is
irreducible.
\end{proof}

These considerations allow us to introduce the dual complex
$\dc(Z'')$ for the exceptional divisor $Z''$ of the morphism
$\alpha\colon Y''\to X$. Since all the strata $Z_\tau$, $\dim
Z_\tau\geq 2$, have irreducible intersections with $X$, we can
formulate the following receipt for finding $\Gamma(Z'')$. First
construct the dual complex $\dc(\Bar{Z''})$, where $\Bar{Z''_i}$ are
the exceptional divisors of the morphism 
$$\Bar{\alpha}\colon V(\Sigma'')\to\C^{n+1}\,,$$ 
$Z''_i=\Bar{Z''_i}\cap Y''$. It can be built
using a hyperplane section $H\cap\Sigma''$ as
described above. Then remove from $\dc(\Bar{Z''})$ all the interiors
of $n$-dimensional cells of $\Sigma''\cap H$ and glue additional
$(n-1)$-simplexes for those 1-dimensional strata $Z_\tau$ which intersect
$Y''$ in more than one point. To make this construction clearer let us
illustrate it by an example (taken from~\cite{Gordon1}).
\begin{example}\label{Ex:Malgrange}
Consider
$$(X,0)=(\{x^8+y^8+z^8+x^2y^2z^2=0\},0)\subset(\C^3,0)\,.$$
This is a non-degenerate isolated surface singularity. Its Newton
diagram is shown in Fig.~\ref{F:Malgrange}, {\it a}.
\begin{center}
\begin{figure}[h]
\begin{picture}(10,6)
\put(1,1){\vector(1,0){4}} \put(4.6,1.15){$\deg(x)$}
\put(1,1){\vector(2,3){2}} \put(3.1,3.9){$\deg(y)$}
\put(1,1){\vector(0,1){4}} \put(1.1,4.9){$\deg(z)$}
\put(1.75,1.75){\line(-1,3){0.75}} \put(1.75,1.75){\line(1,3){0.25}}
\put(1.75,1.75){\line(3,-1){2.25}}
\put(4,1){\line(-4,3){2}} \put(4,1){\line(-1,1){3}}
\put(1,4){\line(2,-3){1}}
\put(4,0.6){$8$} \put(1.9,2.65){$8$} \put(0.7,3.8){$8$}
\put(1.4,1.3){$(2,2,2)$}
\put(2.5,0.4){\it a}
\put(6,1){\line(1,0){3.5}} \put(6,1){\line(2,1){1}}
\put(6,1){\line(1,2){1.75}}
\put(9.5,1){\line(-1,2){1.75}}
\put(7,1.5){\line(1,0){1.5}} \put(7,1.5){\line(1,2){0.75}}
\put(8.5,1.5){\line(2,-1){1}} \put(8.5,1.5){\line(-1,2){0.75}}
\put(7.75,3){\line(0,1){1.5}}
\put(7,1.5){\circle*{0.1}} \put(8.5,1.5){\circle*{0.1}}
\put(7.75,3){\circle*{0.1}}
\put(5.7,0,6){$(1,0,0)$} \put(9.2,0.6){$(0,1,0)$}
\put(7.4,4.6){$(0,0,1)$}
\put(6.8,1.1){$(2,1,1)$} \put(8.5,1.65){$(1,2,1)$}
\put(7.9,2.8){$(1,1,2)$}
\put(7.7,0.4){\it b}
\end{picture}
\caption{Newton diagram for $x^8+y^8+z^8+x^2y^2z^2$.}\label{F:Malgrange}
\end{figure}
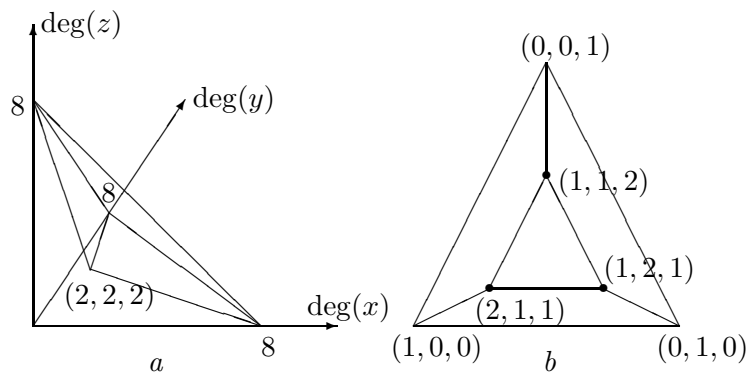
\end{center}
The section $\Sigma'\cap H$ of the first Varchenko subdivision with an
appropriate hyperplane is 
shown in Fig.~\ref{F:Malgrange}, {\it b}. $(1,0,0)$ etc. are primitive
vectors along the corresponding rays; bold points indicate the
exceptional divisors. Three cells of $\Sigma'\cap H$ are not
triangles. We can triangulate them by introducing some vectors into
their interiors (i.~e., by additional toric blow-ups); but by
Proposition~\ref{P:intersection} the corresponding exceptional
divisors does not intersect the strict transform of $X$, thus they will
be deleted in the next step. The complex $\dc(\Bar{Z''})$ is shown in
Fig.~\ref{F:dcMalgrange}, {\it a}. It is just the inner triangle of
$\Sigma'\cap H$.
\begin{center}
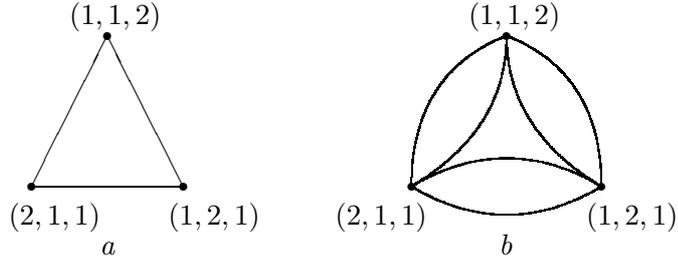
\begin{figure}[h]
\begin{picture}(10,4)
\put(1.5,1){\line(1,0){2}} \put(1.5,1){\line(1,2){1}}
\put(3.5,1){\line(-1,2){1}}
\put(1.5,1){\circle*{0.1}} \put(3.5,1){\circle*{0.1}}
\put(2.5,3){\circle*{0.1}}
\put(1.2,0.5){$(2,1,1)$} \put(3.3,0.5){$(1,2,1)$}
\put(2,3.15){$(1,1,2)$}
\put(2.4,0.1){\it a}
\put(6.5,1){\circle*{0.1}} \put(9,1){\circle*{0.1}}
\put(7.75,3){\circle*{0.1}}
\qbezier(6.5,1)(6.5,2.5)(7.75,3) \qbezier(6.5,1)(7.75,1.75)(7.75,3)
\qbezier(6.5,1)(7.75,0.25)(9,1) \qbezier(6.5,1)(7.75,1.75)(9,1)
\qbezier(9,1)(9,2.5)(7.75,3) \qbezier(9,1)(7.75,1.75)(7.75,3)
\put(5.5,0.5){$(2,1,1)$} \put(8.8,0.5){$(1,2,1)$}
\put(7.25,3.15){$(1,1,2)$}
\put(7.65,0.1){\it b}
\end{picture}
\caption{Dual complex of $x^8+y^8+z^8+x^2y^2z^2$.}\label{F:dcMalgrange}
\end{figure}
\end{center}
To get the dual complex $\dc(Z'')$ we must delete the interior of the
triangle in Fig.~\ref{F:dcMalgrange}, and the only remaining question 
is whether its edges are multiple. In order to find this out let us
consider the birational morphism $V(\Sigma'')\to\C^3$ restricted to
the affine piece $V(\sigma)=\C^3/G$, $\sigma=\langle(2,1,1), (1,2,1),
(1,1,2)\rangle$, and $G$ is a group of order
$$
\begin{vmatrix}
2 & 1 & 1 \\
1 & 2 & 1 \\
1 & 1 & 2
\end{vmatrix}
 =4\,.
$$
Formulae~\eqref{E:tmap} take the form
$$x=x_{1}^{2}y_1z_1,\; y=x_1y_{1}^{2}z_1,\;z=x_1y_1z_{1}^{2}\,.$$
The strict transform $Y''$ of $X$ is
$$x_{1}^{8}+y_{1}^{8}+z_{1}^{8}+1=0\,.$$
The intersection of $\widetilde{Y}''$ with, say, the stratum 
$Z_\tau$: $y_1=z_1=0$
consists of 8 points; but we must account also the group
action. Verification shows that $G=\Z_4$ acting via
$$x_1\to\varepsilon x_1\,,\; y_1\to\varepsilon y_1\,,\; z_1\to\varepsilon
z_1\,,$$
$\varepsilon^4=1$, thus $\sharp (Y''\cap Z_\tau)=2$. Therefore,
the edges of triangle are double. The complex $\dc(Z'')$ is shown in
Fig.~\ref{F:dcMalgrange}, {\it b}. It is homotopy to a bouquet of 4
circles; in particular, rank of $H_1(\dc(Z''))$ is $4$.
\end{example}

Note that constructing $\Sigma''$ can be much easier than finding
$\Sigma$ that gives the embedded resolution. However, knowing
$\Sigma''$ is enough to determine the homotopy type of $\dc(Z)$.
\begin{proposition}\label{P:partres}
Let the singularity $(X,o)$ be such that the property (R) is
satisfied, and let $\Sigma''$ be as described above. Then $\dc(Z'')$
is homotopy equivalent to the dual complex $\dc(Z)$ associated to any
good resolution $(Y,Z)\to(X,o)$.
\end{proposition}
\begin{proof}
The proposition is an easy consequence of the theory of toroidal
embeddings. Indeed, let $Y''$ be the strict transform of $X$ in
$V(\Sigma'')$. Since all cones of $\Sigma'$ are simplicial, $Y''$ has
only toroidal quotient singularities. Then
$Y''\setminus Z''\subset Y''$ is a toroidal embedding. In a similar
way to the proof of Theorem~\ref{T:terminal} we can construct a good
toroidal resolution $Y$ of $Y''$ such that its dual complex $\dc(Z)$
is a subdivision of $\dc(Z'')$. Hence $\dc(Z)$ is homeomorphic to
$\dc(Z'')$. By Theorem~\ref{T:invariant} dual complex of any good
resolution of $X$ is homotopy to $\dc(Z)$ and thus to $\dc(Z'')$.
\end{proof}

Applying Proposition~\ref{P:partres} to Example~\ref{Ex:Malgrange}, we
deduce that $\dc(Z)$ is homotopy equivalent to the bouquet of $4$
circles for any good resolution $Y\supset Z$. This agrees with
\cite{Gordon1} where $(X,0)$ is resolved by a sequence of appropriate
blow-ups. One more application is the following
\begin{corollary}\label{C:Brieskorn}
Let $(X,0)$ be the Brieskorn singularity
$$x_{0}^{a_0}+x_{1}^{a_1}+\dots+x_{n}^{a_n}=0\,.$$
Let $(Y,Z)$ be its good resolution. Then $\dc(Z)$ has trivial homotopy
type.
\end{corollary}
\begin{proof} In view of Proposition~\ref{P:partres} the task becomes
trivial. The first Varchen\-ko subdivision in this case is just the
subdivision of the positive cone $\sigma_+$ by the ray $\langle
w\rangle$, 
$$w=\left(\frac{m}{a_0},\frac{m}{a_1},\dots,\frac{m}{a_n}\right)\,,$$
where $m$ is the least common multiple of the integers $a_0$,
$a_1,\dots,a_n$. Thus $\Sigma'$ consists of simplicial cones and
possesses the property (R). We can put $\Sigma''=\Sigma'$. There is
only one prime exceptional divisor, hence $\dc(Z'')$ is a
point. Therefore $\dc(Z)$ is homotopy equivalent to a point for any
good resolution.
\end{proof}

We have seen that for a partial resolution $Y''\subset V(\Sigma'')$ of
a non-degenerate singularity the complex $\dc(Z'')$ is just the
$(n-1)$-skeleton of $\dc(\Bar{Z''})$ to which, maybe, several
additional $(n-1)$-simplexes have been glued. But $\dc(\Bar{Z''})$ is
homotopy trivial. Indeed,  we can
make all its cones simplicial by introducing new rays only in its
$n$-skeleton. Also, since $\Sigma'$ satisfies the property (R), we can
put all these rays to the interior of $\sigma_+$. Then 
$\dc(\Bar{Z''})$ is obtained from
$\dc(T+\Bar{Z''})$ by deleting the vertices corresponding to the
divisor $T$. This does not change the homotopy type of the dual
complex and, on the other hand, $\dc(T+\Bar{Z''})$ is obviously
homotopy trivial. We come to the following result about which we have
a feeling that it must hold in a much wider situation.
\begin{corollary}\label{C:nondeg}
Let $(X,0)$ be a non-degenerate isolated $n$-dimensional hypersurface
singularity satisfying the property (R). Let $(Y,Z)\to(X,0)$ be a good
resolution. Then all the intermediate homologies of $\dc(Z)$ vanish:
$$H_k(\dc(Z),\Z)=0 \text{ for } 0<k<n-1\,.$$
\end{corollary}


\begin{thebibliography}{24}
\bibitem{AKMW} Abramovich, D., Karu, K., Matsuki, K., 
  W{\l}odarczyk, J. \emph{Torification and factorization of birational 
  maps}, J. Amer. Math. Soc. \textbf{15} (2002), 531--572. 
\bibitem{ACampo} A'Campo, N. \emph{La fonction zeta d'une 
   monodromie}, Comment. Math. Helvetici \textbf{50} (1975), 233--248.
\bibitem{Artin} Artin, M. \emph{On isolated rational singularities
  of surfaces}, Amer. J. Math. \textbf{88} (1966), 129--136.
\bibitem{Danilov} Danilov, V.\,I. \emph{The geometry of toric
  varieties} (Russian), Uspekhi Matem. Nauk \textbf{33}(2) (1978), 85--134.
\bibitem{Elkik} Elkik, R. \emph{Rationalit{\'e} des singularit{\'e}s 
  canoniques}, Invent. Math. \textbf{64} (1981), 1--6.
\bibitem{Fulton} Fulton, W. \emph{Introduction to Toric Varieties},
  Princeton University Press, Princeton, 1993.
\bibitem{Hauser} Encinas, S., Hauser, H. \emph{Strong resolution of
  singularities in characteristic zero},
  Comment. Math. Helv. \textbf{77}(4) (2002), 821--845.
\bibitem{Gordon1} Gordon, G.\,L. \emph{On a simplicial complex 
  associated to the monodromy}, Transactions of the AMS \textbf{261} 
  (1980), 93--101.
\bibitem{Gordon2} Gordon, G.\,L. \emph{On the degeneracy of a spectral
  sequence associated to normal crossings}, Pacific
  J. Math. \textbf{90}(2) (1980), 389--396.
\bibitem{Grauert} Grauert, H. et al. \emph{Several complex 
  variables VII. Sheaf-theoretical methods in complex analysis}, 
  Encyclopedia of Mathematical Sciences. 74. 
  Berlin: Springer-Verlag, 1994.
\bibitem{Grothendieck} Grothendieck, A. \emph{{\'E}l{\'e}ments de
  g{\'e}om{\'e}trie alg{\'e}brique}, Ch. III, Institut des Hautes
  {\'E}tudes Scienttifiques, Publications Math{\'e}matiques No. 11
  (1961).
\bibitem{Hironaka} Hironaka, H. \emph{Resolution of singularities
  of an algebraic variety over a field of characteristic zero},
  Annals of Math. \textbf{79} (1964), 109--326.
\bibitem{Mumford} Kempf, G., Knudsen, F., Mumford, D., Saint-Donat,
  B. \emph{Toroidal embeddings I}, Springer LNM 339, 1973.
\bibitem{Kulikov} Kulikov, Vik.\,S., Kurchanov, P.\,F. \emph{Complex
  algebraic varieties: periods of integrals and Hodge structures},
  Algebraic geometry, III, 1--217. Encyclopaedia Math. Sci., 36,
  Springer, Berlin, 1998.
\bibitem{Matsuki} Matsuki, K. \emph{Lectures on factorization of
  birational maps}, e-preprint: arXiv:math.AG/0002084.
\bibitem{Milnor} Milnor, J. \emph{Singular points of complex
  hypersurfaces}, Annals of Mathematics Studies 61, Princeton
  University Press, Princeton, 1968.
\bibitem{Mori} Mori, S. \emph{On 3-dimensional terminal
  singularities}, Nagoya Math. J. \textbf{98} (1985), 43--66.
\bibitem{RY} Reid, M. \emph{Young person's guide to canonical
  singularities}, Proc. Sympos. Pure Math., 46, Part 1,
  Amer. Math. Soc., Providence, 1987.  
\bibitem{Shokurov} Shokurov, V.\,V. \emph{Complements on surfaces},
  J. Math. Sci. (New York) \textbf{102}(2) (2000), 3876--3932.
\bibitem{Stepanov1} Stepanov, D.\,A. \emph{A note on the dual complex
  associated to a resolution of singularities}, Uspekhi Matem. Nauk
  \textbf{61}(1) (2006), 185--186.
e-preprint: arXiv:math.AG/0509588.
\bibitem{Stepanov2} Stepanov, D.\,A. \emph{A note on resolution of
  rational and hypersurface singularities}, submitted to Proc. of the
  AMS. e-preprint: arXiv:math.AG/0602080.
\bibitem{Thuillier} Thuillier, A. \emph{G{\'e}om{\'e}trie
  toro{\"{\i}}dale et g{\'e}om{\'e}trie analytique non
  Archim{\'e}dienne. Application au type d'homotopie de certains
  sch{\'e}mas formels}, Preprint Nr. 10 (2006), Universit{\"a}t
  Regensburg, Mathematik.
\bibitem{tomDieck} tom Dieck, T. \emph{Transformation groups}, De
  Gruyter studies in mathematics 8. Berlin; New York: de Gruyter, 1987.
\bibitem{Varchenko} Varchenko, A.\,N. \emph{Zeta-function of monodromy
  and Newton's diagram}, Invent. Math. \textbf{37}(3) (1976),
  253--262.
\end{thebibliography}
\end{document}